\documentclass[10pt]{amsart}

\usepackage{amssymb}    % use all AMS fonts
\usepackage{amsmath}    % include some AMS-LaTeX functionality
\usepackage{amsthm}     % AMS style theorem and proof environments
\usepackage{mathrsfs}   % special calligraphic characters like \mathscr{S}
\usepackage{euscript}   % euscript package
\usepackage{epsfig}             % epsfig package, encap postscript
\usepackage[all]{xy}    % Xy-Pic for diagrams

\begin{document}

\newcommand{\minus}{\setminus}
\newcommand{\cross}{\times}
\newcommand{\st}{\mid}
\newcommand{\comp}{\circ}
\newcommand{\norm}{\Arrowvert}
\newcommand{\ra}{\rightarrow}
\newcommand{\R}{\mathbb{R}}
\newcommand{\Q}{\mathbb{Q}}
\newcommand{\Hy}{\mathbb{H}}
\renewcommand{\epsilon}{\varepsilon}
\newcommand{\nin}{\not\in}
\renewcommand{\qed}{\square}
\renewcommand{\phi}{\varphi}
\newcommand{\ri}{\underline{\rho}}
\newcommand{\rs}{\overline{\rho}}
\newcommand{\te}{t \epsilon}
\newcommand{\cc}{d_{cc}}
\newcommand{\fcc}{f_{cc}}
\newcommand{\db}{\overline{d}}
\newcommand{\pt}{\|_t}
\newcommand{\ghlim}{\text{GH-lim}}
\newcommand{\wlim}{\text{$\omega$-lim}}
\newcommand{\Ha}{\mathscr{H}}
\newcommand{\ip}{<\cdot,\cdot>}
\newcommand{\pf}{\noindent {\em Proof: }}
\newcommand{\df}{\underset{def}{\equiv}} 
\newtheorem{Pro}{Proposition}[section]
\newtheorem{Lem}[Pro]{Lemma}
\newtheorem{Sub}[Pro]{Sublemma}
\newtheorem{Thm}[Pro]{Theorem}
\newtheorem{MThm}{Theorem}
\renewcommand{\theMThm}{\Alph{MThm}}
\newtheorem{Rem}{Remark}
\newtheorem{Def}[Pro]{Definition}
\newtheorem{Not}{Notation}
\newtheorem{Cor}[Pro]{Corollary}
\newtheorem{Ithm}{Theorem}
\renewcommand{\theIthm}{\arabic{chapter}.\arabic{Ithm}}
\newtheorem{Idef}{Definition}
\renewcommand{\theIdef}{\arabic{chapter}.\arabic{Idef}}

\title{The large scale geometry of Nilpotent Lie groups}
\author{Scott D. Pauls}
\address{Rice University, Houston, TX, 77005}
\email{pauls@math.rice.edu}

\begin{abstract}
In this paper, we prove results concerning the large scale geometry of
connected, simply connected nilpotent Lie groups equipped with left
invariant Riemannian metrics.  Precisely, we prove that there do not
exist quasi-isometric embeddings of such a nilpotent Lie group into
either a $CAT_0$ metric space or an Alexandrov metric space with
curvature bounded below.  The main
technical aspect of this  
work is the proof of a limited metric differentiability of Lipschitz maps
between connected graded nilpotent Lie groups equipped with left
invariant Carnot-Carath\'eodory metrics and complete metric spaces.   
\end{abstract}

\maketitle
\section{Introduction}
In this paper we investigate the large scale geometry of connected
nilpotent Lie groups equipped with left invariant Riemannian metrics
by studying their quasi-isometric embeddings into various metric
spaces.  Let $(N,g)$ be a connected nilpotent Lie group
with a left invariant Riemannian metric and $d$ be the induced
distance function on $N$.  If $(X,d_X)$ is a complete metric space,
then $f:N \ra X$ is an $(L,C)$-quasi-isometric embedding if, for all $x,y
\in N$,
\[\frac{1}{L}d(x,y) -C \le d_X(f(x),f(y)) \le Ld(x,y) +C \]
After studying certain invariants of these maps, we prove two main
applications:

\begin{MThm}\label{noCAT0}  There do not exist quasi-isometric
embeddings of a connected nonabelian nilpotent Lie group equipped
with a left invariant Riemannian metric into a $CAT_0$ metric space. 
\end{MThm}

\begin{MThm}\label{noCBB0}  There do not exist quasi-isometric
embeddings of a connected nonabelian nilpotent Lie group equipped
with a left invariant Riemannian metric into a $CBB_0$ metric space. 
\end{MThm}
$CAT_\kappa$ (resp. $CBB_\kappa$) metric spaces are spaces of
curvature bounded above (resp. below) by $\kappa$ is the sense of
Topanogov triangle comparison.   Thus, a $CAT_0$ metric space is a
generalized space of nonpositive 
curvature while a $CBB_0$ metric space is a generalized space of
nonnegative curvature.  These include simply connected Riemannian
manifolds of 
nonpositive and nonnegative curvature respectively.  In the
literature, $CAT_0$ spaces are also called Hadamard spaces and
$CBB_\kappa$ spaces 
are called  
Alexandrov spaces (with curvature bounded below).  We retain the earlier notion for simplicity and
consistency.  We will discuss
these in detail later in the 
paper.  In \cite{Wolf:nil}, Wolf showed that a
connected nonabelian nilpotent Lie group equipped with a left
invariant Riemannian metric must contain 2-planes of both positive and
negative curvatures.  Theorems \ref{noCAT0} and \ref{noCBB0} can be
viewed as a large 
scale analogue of this theorem of J. Wolf.  

While the large scale characteristics of such nilpotent Lie groups
may be interesting in their own right, we shall see that such
investigations reduce to problems which are motivated by the now
standard arguments used in the proof of Mostow's rigidity theorem and
its many extensions (see \cite{Mostow}, \cite{Ballm:rigid},
\cite{Burns:Spatzier}, \cite{Pansu},
\cite{Hamen:rank1} as well as many others).  In the proofs of these
results, one attempts to show that two candidate spaces (e.g. two
compact constant negative curvature spaces with isomorphic
fundamental groups) are isometric by exhibiting an equivariant
quasi-isometry between 
their universal covers and showing that the existence of such a quasi-isometry
must imply the existence of a true equivariant isometry.  A
common element in many such proofs is to reduce to an examination of
the ideal boundaries of the two spaces in question.  The
quasi-isometries between the two spaces induce quasiconformal maps on the boundaries
which then become the focus of the study.  In the cases of the papers
mentioned above, the tangent cones at points on the ideal boundaries
are isometric to connected simply connected graded nilpotent Lie group
equipped with left invariant 
Carnot-Carath\'eodory metrics.  Thus, a local analysis of the ideal
boundaries involves an examination of the geometry of graded nilpotent
Lie groups with left invariant Carnot-Carath\'eodory metrics.  This
provides a natural motivation for considering the geometry of such
spaces and, in particular, of studying the quasi-conformal maps
between them.  It is a natural generalization to instead consider the
quasi-conformal embeddings of such spaces.  Unfortunately, the study
of such embeddings 
seems quite intractable but does provide motivation for various lines of
study - for example the study of biLipschitz embeddings of such spaces. 
As we shall see, the proofs of theorems \ref{noCAT0} and \ref{noCBB0}
rest on a local analysis of biLipschitz embeddings of graded nilpotent Lie
groups with left invariant Carnot-Carath\'eodory metrics into
complete metric spaces which are {\it locally} $CAT_\kappa$ or
$CBB_\kappa$.  For metric spaces which are locally $CAT_\kappa$,
i.e. for each point there is a closed ball about that point which is
itself $CAT_\kappa$, we use the notation $CBA_\kappa$ (for ``curvature
bounded above'').  Note that, for a lower curvature bound, these
distinctions are not necessary.  In particular, we prove an intermediate
theorem from which theorems \ref{noCAT0} and \ref{noCBB0} follow:
\begin{MThm}\label{Inter}  Let $G$ be a connected simply connected
  graded nilpotent Lie group equipped with a left invariant
  Carnot-Carath\`eodory metric and $U \subset G$ be an open set.  Then
  $U$ does not admit a biLipschitz embedding into any $CBA_\kappa$
  metric space or into any $CBB_\kappa$ metric space.
\end{MThm}  
To prove the main theorems from this intermediate theorem, we study
some biLipschitz embedded invariants of graded 
nilpotent Lie groups with left invariant 
Carnot-Carath\'eodory metrics and use them to build local obstructions
to quasi-isometric embeddings of nonabelian nilpotent Lie groups with left
invariant Riemannian metrics into various metric spaces.  Next we give
a brief outline of the argument and the paper and show how theorems
\ref{noCAT0} and \ref{noCBB0} follow from theorem \ref{Inter}.  

Let $N$ be a connected simply connected nilpotent Lie group with a
left invariant 
Riemannian metric $g$ and let $d$ be the induced distance function on
$N$.  If $f:N \ra X$ is some $(L,C)$-quasi-isometric embedding of $N$ into a
complete metric space $X$, we first consider asymptotic cones of
$N$ and $X$ and a map between them, denoted $F$, derived from $f$
and the coning procedure.  The asymptotic cone of a metric space
$(X,d_X)$ is a limit metric space of the pointed dilated spaces
$(X,x_i,\lambda_i d_X)$ for some sequence $\lambda_i \ra 0$.  Since
such spaces do not necessarily converge in the Gromov-Hausdorff
topology, we use Gromov's ultrafilter construction, choosing a
nonprincipal ultrafilter $\omega$ (see section \ref{Tangcones} for a
definition) to form the asymptotic cone $(C_\infty^\omega
X,d_X^\omega)$.  For precise details of the construction, we refer
the reader to either \cite{Gromov:AIT} or \cite{KleinerLeeb}.  For our
purposes, there are a few keys pieces of information.  First, if
$(N,g)$ is a nilpotent Lie group with a left invariant Riemannian
metric, then Pansu (\cite{Pansu}) proved that the asymptotic cone is
unique and isometric to $(\widetilde{N},\cc)$, a graded nilpotent Lie
group equipped with a left invariant Carnot-Carath\'eodory metric (see
section \ref{CC} for definitions).  Second, an asymptotic cone of a
$CAT_0$ space is also $CAT_0$ and an asymptotic cone to a $CBB_0$
space is also $CBB_0$.  These facts (and references to their proofs)
are reviewed in section \ref{Curv}.  Third, the quasi-isometric
embedding $f$ gives rise to an L-biLipschitz map $F$ between the cones
$(\widetilde{N},\cc)$ and $(C_\infty^\omega X, d_X^\omega)$.  The
construction of this map depends on the choice of ultrafilter.  Thus,
to prove theorems A and B, we must prove instead theorem \ref{Inter}.
The paper is 
devoted to proving this intermediate theorem.

The main
technical result in the paper is 
that, in an appropriate sense, such an L-biLipschitz map is differentiable
in certain directions 
almost everywhere.  This is a generalization of a classical theorem
of Rademacher. Our proof follows the same lines as a covering argument
used in  
Kleiner's proof of the differentiability of Lipschitz maps from $\R^n$
into metric spaces (in \cite{Kleiner2}), which is a special case of
differentiability theory of section 1.9 in \cite{Korevaar:Schoen}.  In
\cite{Kir}, Kirchheim shows the metric
differentiability of Lipschitz maps $f:\R^n \ra X$ where the targets
are complete metric spaces.  Our differentiability result is an
extension of Kirchheim's work.  The interested reader should also
consult \cite{Pansu} and \cite{Mostow:Margulis} for results concerning
the differentiability of quasiconformal maps between
Carnot-Carath\'eodory spaces.

To prove theorem \ref{Inter}, we construct a tangent cone of $F$ at a
point of differentiability. Using this tangent cone,
we can compare the local geometry of the two asymptotic cones.  The
local geometry of connected graded nilpotent Lie groups with left 
invariant Carnot-Carath\'eodory metrics is well understood and the
differentiability allows us to ``push forward'' this structure to a
tangent cone to the asymptotic cone of $X$.  When this asymptotic cone
has additional structure, this allows us to measure the compatibility
of these two objects.  In the case of Theorem \ref{Inter}, this
structure provides estimates on the rate of 
growth of the spread between two tangent cone geodesics which show
that they spread 
apart sublinearly but do not remain a bounded distance from one
another.  Comparing this to the spread of geodesics in the tangent
cone to either a $CAT_0$ or $CBB_0$ space, we derive a contradiction
which proves Theorem \ref{Inter}.

Sections \ref{Nilrev} and \ref{CC} are reviews 
of the constructions mentioned above and discuss, respectively,
nilpotent Lie groups and  Carnot-Carath\'eodory metrics.  Section 
\ref{secdiff} contains the proof of the limited differentiability of biLipschitz
embeddings of connected graded nilpotent Lie groups with left
invariant Carnot-Carath\'eodory metrics into complete metrics
spaces.  The main goal of this section is to prove theorem
\ref{KScc}.  Section \ref{Tangcones} reviews the tangent cone
construction and interprets theorem \ref{KScc} as a statement about
maps between tangent cones.  Section \ref{Curv} reviews the
definition and some properties of $CAT_0$ and $CBB_0$ spaces.  In section
\ref{end} we prove Theorem C, that there do not exist
biLipschitz embeddings of a nonabelian connected graded nilpotent Lie
group equipped with a Carnot-Carath\'eodory metric into either a
$CBA_\kappa$ or $CBB_\kappa$ space. 

The author wishes to thank Bruce Kleiner for suggesting this line of
work as well as for many helpful discussions.  Thanks are also due to
Chris Croke for helpful discussions of this work and for many hours of
help with the preparation of this document.

\section{Nilpotent Lie Groups}\label{Nilrev}
In this section, we recall some basic definitions of nilpotent Lie
groups and some associated structures.  Let $G$ be a connected Lie group and
$\mathfrak{g}$ be its Lie algebra equipped with $[\cdot,\cdot]$ a
bracket.  For the rest of this exposition, we assume that all Lie
groups are connected and simply connected.
\begin{Def}  Let $\mathfrak{g}$ be a Lie algebra over a field $K$, we
define the {  descending central sequence}, $\{C^k(\mathfrak{g})\}$
\index{$C^k(\mathfrak{g})$} as follows:
\begin{equation*}
\begin{split}
C^0(\mathfrak{g}) &= \mathfrak{g} \\
C^{k+1}(\mathfrak{g}) &= [C^k(\mathfrak{g}),\mathfrak{g}]
\end{split}
\end{equation*}
\end{Def}
We use the descending central sequence to define nilpotency:
\begin{Def}A Lie group $G$ is called
{  nilpotent} if there exists an integer $n$ such that
$C^{n}(\mathfrak{g})=0$.  The smallest such $n$ is called the degree of
nilpotency of $G$.  
\end{Def}
Recall that the exponential map, $v \mapsto e^{v}$, is a
diffeomorphism for simply connected nilpotent Lie groups.  For this
paper, we will use the following notation:  if $e^v$ and $e^w$ are
exponential images of the Lie algebra elements $v$ and $w$, then th e
product, $e^ve^w$ is also the exponential image of a Lie algebra
element which we denote $v \circledcirc w$.  In other words,
$e^ve^w=e^{v \circledcirc w}$. 

For this work, Pansu's theorem tells us that after taking asymptotic
cones, graded nilpotent Lie groups are our main objects of study.

\begin{Def} \index{graded Lie algebra}  A nilpotent Lie group is
  graded if its Lie algebra comes with a grading.  A grading for a Lie
  algebra $\mathfrak{g}$ is a decomposition 
\[\mathfrak{g} = \oplus_i \mathcal{V}^i\]
where the subspaces $\mathcal{V}^i$ satisfy the condition
\[[\mathcal{V}^i,\mathcal{V}^j] \subset \mathcal{V}^{i+j}\]
for all $i,j \in \mathbb{Z}$.  
\end{Def}
 
\noindent
{\em Example: }The three dimensional Heisenberg group is defined by giving a three
dimensional Lie algebra spanned by vectors $X, Y, Z$ with the only
nontrivial bracket relation given by $[X,Y]=Z$.  
Thus, the Lie algebra has a grading
$\mathfrak{h}=\mathcal{V}^1 \oplus \mathcal{V}^2$ where
$\mathcal{V}^1=span\{X,Y\}$ and $\mathcal{V}^2=span\{Z\}$.  

Every graded nilpotent Lie group possesses a self-similarity.
\begin{Def}\label{htdef}\index{$h_t$}  Let $g\in G$ be an
element in a graded
nilpotent Lie group of nilpotency degree $k+1$ such that
$g=e^{g_1+g_2+...+g_k}$ where $g_i \in 
\mathcal{V}^i$.  For $t \ge 0$, consider the $1$-parameter group of
automorphisms, $\{h_t\}$, of $G$ given by  
\[ h_t(g)=e^{tg_1+t^2g_2+...+t^kg_k}\]
We may sometimes analyze this action of Lie algebra elements.  Thus we
define $dh_t(v_1+v_2+ ... +v_k)=tv_1+t^2v_2+...+t^kv_k$.  If we wish
to consider $h_{-|t|}g$ we use the convention that $ h_{-|t|}g=h_{|t|}g^{-1}
$ (Please note that this convention is not standard in the literature).
\end{Def}

For general Lie groups, the Campbell-Baker-Hausdorff formula provides
a method for computing the
product of two nearby group elements.  We state it in the special case when
$N$ is a simply connected graded nilpotent Lie group.  For a more
general version and, for  
that matter, a proof, see \cite{Vara} section 2.15.

\begin{Thm}[Campbell-Baker-Hausdorff Formula]\label{CBH}  Let $N$ be a
simply connected graded nilpotent 
Lie group of nilpotency degree $k+1$, let $\mathfrak{n}$ 
be its Lie algebra and let $e^\cdot:\mathfrak{n} \ra N$ denote the
exponential map at the identity.  Then, given $X,Y \in \mathfrak{n}$
\begin{equation*}
X \circledcirc Y =\sum_{i=1}^{k-1} \frac{(-1)^{i+1}}{i} 
\sum_{\begin{substack}
                1 \le j \le i, \newline \newline
                p_j+q_j >0
                \end{substack}} 
\frac{\sum_{j=1}^i
(p_j+q_j)^{-1}}{p_1!q_1!...p_n!q_n!} (ad X)^{p_1}(ad Y)^{q_1} ... (ad
X)^{p_n}(ad Y)^{q_n-1}Y
\end{equation*}
where $(ad X)Y=[X,Y]$ and if $q_n=0$ then the last term in the sum is 
$(ad X)^{p_1}(ad Y)^{q_1} ... (ad X)^{p_n-1}X$.  In other words, we
have the formula:
\[e^Xe^Y=e^{X \circledcirc Y}=e^{X+Y+\frac{1}{2}[X,Y]+\frac{1}{12}[X,[X,Y]]-
\frac{1}{12}[Y,[X,Y]]+ ...}\]
\end{Thm}\index{Campbell-Baker-Hausdorff formula}

In the case that $N$ is a graded nilpotent Lie group, we use the
following convention. Let $C_i(X,Y)$ be the pieces of $X \circledcirc
Y$ which lie in $\mathcal{V}^i$.  Thus, if $X,Y \in \mathcal{V}^1$,
$C_1(X,Y)=X+Y$, $C_2(X,y)=\frac{1}{2}[X,Y]$ and so on.  Thus, the
Campbell-Baker-Hausdorff formula yields $e^{X \circledcirc
Y}=e^{C_1(X,Y)+...+C_k(X,Y)}$.  Using the Campbell-Baker-Hausdorff
formula, we easily confirm that $h_t$ is a
group automorphism.  

In subsequent sections, we
will use repeatedly several properties of $h_t$.  We group them here
as a lemma for easy reference:
\begin{Lem}\label{ht}  Let $N$ be a graded nilpotent Lie group and let
$h_t$ be the automorphism described above.
\begin{itemize}
\item For $n \in N$, $h_th_s n=h_{ts}n$.
\item For $v \in \mathcal{V}^1$, $h_te^vh_se^v=h_{t+s}e^v$.
\item If $dL_g$ is the derivative of the left translation map $n
\mapsto gn$, then $dh_t \comp dL_n = dL_{h_tn} \comp dh_t$.
\item $dh_t(v \circledcirc w)=dh_t(v) \circledcirc dh_t(w)$
\end{itemize}
\end{Lem}
\pf The first two are simple exercises using the
Campbell-Baker-Hausdorff formula.  The last two are restatements of the
fact that $h_t$ is a group automorphism in terms of $dh_t$.     
$\qed$

\section{Carnot-Carath\'eodory metrics}\label{CC}
Carnot-Carath\'eodory metrics arise in a variety of settings and have
been the object of a great deal of study.  In this work, we examine
graded nilpotent Lie groups equipped with Carnot-Carath\'eodory metrics and
hence the definitions and propositions will be tailored for this small
set of spaces.  For a more general and comprehensive introduction to
the theory, the interested reader is encouraged to consult
\cite{Strichartz}, \cite{Gromov:CC}, and  \cite{CCgeometry:Risler}.
The rough idea of a Carnot-Carath\'eodory metric is to measure
the distance between two points by taking the infimum of the length of
piecewise $C^1$ paths connecting the points over a restricted set of
paths, those 
tangent to some subbundle of the tangent bundle.  More precisely (in
our setting), if $N$ is a connected graded nilpotent Lie group with Lie algebra
grading $\mathcal{V}^1 \oplus ... \oplus \mathcal{V}^k$, we define a
left invariant subbundle of the tangent bundle, $\mathcal{V}$, by
letting the fiber at a point be the left
translate of the subspace $\mathcal{V}^1$.  In the rest of the paper,
we will, by abuse of notation, refer to the subbundle and the fiber by
$\mathcal{V}$.  If we place a norm $|\cdot|$ on each fiber of
$\mathcal{V}$, we define the Carnot-Carath\'eodory distance between
$x,y \in N$ by 
\[\cc(x,y) = \inf_\gamma \biggl \{ \int |\gamma'(t)|dt \bigg | \gamma \text{
connects $x$ to $y$ and $\gamma' \subset \mathcal{V}$} \biggr \}\]
If no such paths exist, we set $\cc(x,y)=\infty$.  We call a path
whose tangent is in the subbundle $\mathcal{V}$ at each point a
horizontal path.  Note that $\cc$ defines a left
invariant metric on 
$N$.  In general, if $| \cdot |$ is merely a norm, we call the
quadruple $(N,\mathcal{V},|\cdot|,\cc)$ a Carnot-Carath\'eodory group
or a sub-Finsler group (if we wish to emphasize that $|\cdot|$ is only
a norm) and $\cc$ a Carnot-Carath\'eodory or sub-Finsler metric.  If
$|\cdot|$ arises from an inner product, $\ip$, on $\mathcal{V}$, we
denote the resulting space by $(N,\mathcal{V},\ip,\cc)$ and we call it
a sub-Riemannian group.  This nomenclature and notation is somewhat
different from the literature which often interchanges the terms
Carnot-Carath\'eodory metric, sub-Riemannian metric, and singular
Riemannian metric freely.  In this exposition, any object adorned with
a ``cc'' will be constructed to be with respect to $\cc$.  For example,
$B_{cc}(x,r)$ is the $\cc$-ball of radius $r$ with center $x$.

It is a consequence of Chow's theorem (\cite{Chow} or see
\cite{Gromov:CC} for a exhaustive discussion of this fact) that any
two points in $(N, \mathcal{V},|\cdot|,\cc)$ can be connected with a
smooth horizontal path.  This leaves open the question of geodesics in
Carnot-Carath\'eodory spaces.  We call a segment a geodesic if it is a
local length minimizer.  Again, in the literature, there are many
different terms/meanings for geodesic paths and/or length minimizing
paths depending on whether the geodesics satisfy certain
differential equations (see, for example,  \cite{Suss:Liu}
and \cite{Montgomery}).  For this exposition, these issues are not
relevant and hence we will use the simplified language.  In a
sub-Riemannian group $(N,\mathcal{V},\ip,\cc)$, there is a nice class
of geodesics which will be useful in later sections.
\begin{Def}Let $(N,\mathcal{V},\ip,\cc)$ be a
left invariant sub-Riemannian metric where $N$ is a connected graded
nilpotent Lie group and fix $n \in N$, $v \in \mathcal{V}$.
Let $V$ be a smooth left  
invariant vector field on $N$ with $V(e^0)=v$.  Let $\gamma$
be an integral curve of 
$V$ such that $\gamma(0)=n$.  We call $\gamma$ a {radial geodesic}
eminating from $n$.  If a basepoint is understood, we shall simply
call $\gamma$ a radial geodesic.
\end{Def} 

Note that radial geodesics are actually geodesics of $\cc$.  Indeed,
considering the quotient $N/[N,N]$ endowed with the left invariant
metric induced by the $cc$ inner product on $\mathcal{V}$, the
quotient map $N \ra N/[N,N]$ is $1$-Lipschitz and takes radial
``geodesics'' to geodesics in the quotient and preserves their length.
Thus, these are geodesics of $\cc$ as well.  In addition, we see that
$h_t$ is a homothety of $\cc$.  Since $h_t$ is an automorphism of $N$
and $dh_t$ preserves $\mathcal{V}$ and acts as a homothety on
$\mathcal{V}$, we see that $h_t$ is a homothety on $(N,\cc)$.

In \cite{Mitchell}, Mitchell describes the local structure of more
general Carnot-Carath\'eodory spaces.  If \newline $(N,\mathcal{V},\ip,\cc)$ is
a sub-Riemannian group, the local structure is much easier to
determine given the existence of a left invariant
metric, a homothety of that metric, and the uniqueness of the Haar
measure.  The following results can either be thought of as Mitchell's
results in this special case or as consequences the extra structure in
this special case.

\begin{Pro}\label{Mitchell}  Let $(N,\mathcal{V},\ip,\cc)$
be a sub-Riemannian group.  Then,
\begin{enumerate}
\item For any $n \in N$, the tangent cone at $n$,
$(C_nN,\overline{d}_{cc})$ exists, is unique, and is isometric to
$(N,\cc)$.
\item Let $Q=\sum_{i=1}^k i \dim \mathcal{V}^i$.  Then $Q$ is the
Hausdorff dimension of $(N,\cc)$.
\item Let $\ip_R$ be a Riemannian completion of the inner product
$\ip$ on all of $\mathfrak{n}$ such that the grading is orthogonal.
Let ``dvol'' denote the Riemannian volume measure associated to this
left invariant Riemannian metric.  Then, if $\Ha^Q$ denotes the $Q$
dimensional Hausdorff measure, $\Ha^Q$ and $dvol$ are commensurable.
In other words, there exist constants $c_1,c_2$ such that
\begin{equation}\label{CCvRball}
c_1 r^Q \le \int_{B_{cc}(x,r)} dvol \le c_2 r^Q 
\end{equation}
\end{enumerate}
\end{Pro}

One should note that in Mitchell's work on more general
Carnot-Carath\'eodory spaces, the constants $c_1,c_2$ must be
recalculated depending on a choice of compact set and and
the the Hausdorff and Riemannian measures are only commensurable on
compact sets.  However, we now see that in our
setting, the left invariance of the metric and the presence of a
homothety allows us have global constants for the volume and Hausdorff
measure comparisons.  First, we start with an illustrative example:
\newline
\newline

\noindent
{\em Example:}  Let $H^3$ denote the three dimensional
Heisenberg group with a left-invariant \newline sub-Riemannian metric.  It is 
well known that the Haar measure on $H^3$ is, in fact, the four
dimensional Hausdorff measure with respect to $\cc$ on $H^3$ (see
\cite{Pansu:balls}).  This confirms the formula 
from Mitchell's theorem.  Indeed,
\begin{equation*}
\begin{split}
V^1=\text{span}(\{X,Y\}), 
V^2=\text{span}(\{X,Y,Z\}), 
V^3=\emptyset\\
\text{therefore, }
Q=1(2-0)+2(3-2)=4
\end{split}
\end{equation*}
In addition, we can quickly see the volume estimate inequality
(\ref{CCvRball}) using the homothety $h_r$.  Consider what $h_r$, with
$r>0$, does to 
the Riemannian volume.  Denote by $g$ the Riemannian metric given by
the inner product $\ip_R$ described above.  At the identity, the
metric $h_r^*g$, written in terms of the left invariant 
vector fields $X,Y$ and $Z$ looks like
\begin{gather*}
\begin{pmatrix}
r^2 & 0 & 0 \\
0 & r^2 & 0 \\
0 & 0 & r^4\\
\end{pmatrix}
\end{gather*}
Thus, the Riemannian volume element determined by the metric $g$ is
multiplied by $r^4$ under $h_r^*$.  Thus, picking $C$ to work for a particular value
of $r_0$, we can then use the dilation $h_r$ to scale the picture so
that the same constant works for any $r$ value larger than zero. 

This argument actually works for more general graded nilpotent Lie
groups yielding the same comparative scaling between the
Carnot-Carath\'eodory metric and the Riemannian metric.  

In the proof of the limited differentiability of biLipschitz maps, we
will need to be able to manipulate covers of 
sets in graded nilpotent Lie groups while using the measure $\Ha^Q$.
In particular, we want to be able to take the analogue of a Vitali
cover of a measurable set $S$ and find a countable disjoint subcover
which still covers almost all of $S$ with respect to $\Ha^Q$.  For
this section, we will use the language of \cite{Fed} (section 2.8).
Assume that $N$ is a connected graded nilpotent Lie group with a
sub-Riemannian metric, $\cc$.  Let $\Ha^Q$ be the $Q$
dimensional Hausdorff measure associated to $\cc$ as in last section.
Let $\mathscr{C}$ be the collection of closed cc-balls in $N$.  Then,
mimicing the standard proof in $\R^n$ we prove the following lemma.
\begin{Lem}\label{Vitali}  $\mathscr{C}$ is a Vitali cover.  In other
  words, given a measurable set $A$ covered by a fine subcollection of
  $\mathcal{C}$, we may find a countable dijoint subcover,
  $\mathscr{C}_A$ such that $\Ha^Q_{cc}(A \minus \cup_{C \in
  \mathscr{C}_A} C)=0$.
\end{Lem}

\section{Differentiability of biLipschitz maps}\label{secdiff}
Let $(N,\mathcal{V},\ip,\cc)$ be a graded nilpotent Lie group of
nilpotency degree $k+1$ with a
sub-Riemannian metric.  To set up the notation, we let $F:(N,\cc) \ra
(X,d_X)$ be an 
L-biLipschitz map of $N$ into a complete metric space.  Let $d$ be
$F^*d_X$, the
metric $d_X$ pulled back through $F$.  We want to show that,
appropriately defined, $d$ is differentiable in $\mathcal{V}$
directions at most points.  As discussed in the introduction, to prove
this, we follow the style of the covering argument in the appendix of 
\cite{Kleiner2} used to show a modified version one of Korevaar and
Schoen's results in \cite{Korevaar:Schoen}.  One can also consider
this a version of Kirchheim's metric differentiability in \cite{Kir}
although the argument is different.    
\subsection{Definitions and Preliminary lemmas}
In examining differentiability, we consider an analog of the
directional derivative at $x$ in the direction of $v \in
\mathcal{V}$ 
by looking at difference quotients along the paths $yh_te^v$ for $t>0$:
\[ {d(y,yh_te^v) \over t} \text{ where } y \in B_{cc}(x,t)\]
We consider whether such a ``directional derivative'' exists (i.e. the
limit as $t \ra 0$ exists) by examining the ``lower'' and ``upper''
limits.  If $U$ is some open set in $N$, we define these limits as
follows.  
\begin{Def}\label{lowerd}\index{derivate!lower}\index{$\ri$}
The {  lower derivate} of $d$ is a function $\ri:U \cross \mathfrak{n} \ra
\R$ defined by:
\[\ri(x,v)= {\underline{\lim}}_{t \ra 0^+}\biggl \{ {d(y,yh_te^v) \over t}
\bigg | y \in B_{cc}(x,t)
\biggr \}\]
\end{Def}
\begin{Def}\label{upperd}\index{derivate!upper}\index{$\rs$}
The {  upper derivate} of $d$ is a function $\rs:U \cross \mathfrak{n} \ra
\R$ defined by:
\[\rs(x,v)= {\overline{\lim}}_{t \ra 0^+}\biggl \{ {d(y,yh_te^v) \over t}
\bigg | y \in B_{cc}(x,t)
\biggr \}\]
\end{Def}
\noindent
Note that both $\ri$ and $\rs$ are measurable functions.  Properly
viewed, they are just the $\liminf$ and $\limsup$ of measurable
functions of more than one variable (see \cite{Fed} page 152).  

To understand and estimate the lower and upper derivates, we want to
estimate the quotients for small values of $t$.  To do this, we
construct special boxes that aid in the estimates.  
For each $v \in \mathcal{V}$, denote by $V$ the left invariant vector
field on $N$ determined by $v$.  Choose $\ip_R$ \index{$\ip_R$} to be
a Riemannian completion of the sub-Riemannian metric on
$\mathcal{V}$ such that the grading is orthogonal.
\begin{Def}\label{boxend} Fix $v \in \mathcal{V}$ and let
$V^\perp$ be the orthogonal 
complement of $v$ in $\mathcal{V}$ with respect to
$\ip$.  Let $end(\epsilon,v)=\{w_1 + ... + w_k \in V^\perp \oplus \mathcal{V}^2
\oplus ... \oplus \mathcal{V}^k | w_1 \in V^\perp, w_i \in
\mathcal{V}^i \text{ for $i>1$}, <w_j,w_j>_R^{\frac{1}{2}} <
\epsilon^j, \forall j\}$.  Now, let
\[ End(x,v,\epsilon)=xe^{end(\epsilon,v)}\]
$End(x,v,\epsilon)$ forms the end of our box.
\end{Def} \index{$End(x,v,\epsilon)$}

Note that $e^{end(\epsilon,v)}$ is the exponential image of a
carefully constructed rectangle which should be thought of as a square
of side length $2\epsilon$.  Recalling theorem 3.6 in \cite{Kar:Mont}
or one of the various formulations of the ``ball-box'' estimates (see,
for example \cite{Gromov:CC}),
if $v=v_1+...+v_k \in
\mathfrak{n}$ then $\sum_{i=1}^k <v_i,v_i>_R^{\frac{1}{2i}} \le
A\cc(e^0,e^v)$ where $A$ is a constant depending only on the choice of
$\ip_R$ and the structure constants of the Lie group $N$.  Clearly,
the constant $A$ is invariant under scaling by $h_t$:  if $v'=dh_t v$
and $v_i'$ are the components of $v'$ in the various pieces of the
grading, then $\sum_{i=1}^k < v'_i,v'_i>^\frac{1}{2i} = \sum_{i=1}^kt <
v_i,v_i>^\frac{1}{2i} \le A t \cc(e^0,e^v)=A \cc(e^0,e^{dh_t v})$.  
\begin{Def}\label{box} Define the estimate boxes by
\[Box(x,v,\epsilon) = \bigcup_{z \in End(x,v,\epsilon)} \bigl
(\bigcup_{s \in 
[0,1]} zh_se^v \bigr )\]
\end{Def}\index{$Box(x,v,\epsilon)$}
As in Kleiner's argument, one first thickens $xh_te^v|_{t \in [0,1]}$
into an estimate box, $Box(x,v,\epsilon)$, and then uses smaller
boxes, $Box(y,tv,t\epsilon)$ with $y$ close to $x$, to cover it and
help estimate the 
distance between the ends of the estimate box.  To use this
argument, we must ensure that the boxes we 
have constructed have metric properties (with respect to $\cc$) similar
to the boxes in $\R^n$, at least on the small scale.  One important
technical detail is to control the distances between points in the
``far'' end of the box with respect to corresponding distances in the
``near'' end.  Let $z=yh_{t \epsilon}e^{z_1+z_2+...+z_k} \in
End(y,tv,t\epsilon) 
\subset Box(y,tv,t \epsilon)$ with $\cc(e^0,e^{z_1+z_2+...+z_k}) \le
1$ and $v \in \mathcal{V}$.  Then,
\begin{equation}\label{spread}
\begin{split}
\cc(yh_te^v,zh_te^v) &= \cc(e^0,e^{-tv}e^{\te z_1+ (\te)^2z_2 +
... +(\te)^k z_k}e^{tv})\\
&= \cc(e^0,e^{\te z_1 + ((\te)^2z_2+2t^2 \epsilon c_2[z_1,v])+
... +((\te)^kz_k + 2t^k \epsilon c_{k-1}[z_{k-1},v]+ ...)}) \\  
&= t \cc(e^0,e^{\epsilon z_1 + ((\epsilon)^2z_2+2 \epsilon c_2[z_1,v])+
... +((\epsilon)^kz_k + 2 \epsilon c_{k-1}[z_{k-1},v]+ ...)}) \\
&\le tC(\epsilon)
\end{split}
\end{equation}
Here, $C(\epsilon)$ is a constant that depends on $\epsilon$ and goes
to zero with $\epsilon$.  This constant exists by the compactness of
the closed unit CC-ball.  The first equality is due to the left
invariance of $\cc$.  The second is a use of the
Campbell-Baker-Hausdorff formula.  The third follows from the fact the $h_t$ is a 
homothety of $\cc$.  

As in Kleiner's argument, we will prove a relation between the upper
and lower derivates by covering boxes which thicken the segment
$xh_te^v|_{t \in [0,1]}$ with well controlled smaller boxes.  To do
this, we make the following estimates.  

Let $\mathcal{B}(x,v,\beta,\epsilon)=\{Box(y,tv,t\beta) \bigg |
|{d(y,yh_te^{v}) \over t} - \ri(x,v)| < \epsilon, y \in B_{cc}(x,t)\}$.\index{$\mathcal{B}(x,v,\beta,\epsilon)$}
In other words, members of $\mathcal{B}$ are boxes of some variable
height $t\beta$ 
such that the difference quotients of the endpoints of the middle
segment are infimal up to an error of $\epsilon$.  In fact, we get a
specific estimate on the difference quotient between points on the
ends.
\begin{Lem}\label{endest}  Let $u_1,u_2$ be points on opposite ends of
$Box(y,t{v},t \beta) \in
\mathcal{B}(x,v,\beta,\epsilon)$.  Then, \[{d(u_1,u_2) \over t} \le
C_1(\beta) + \ri(x,v) + \epsilon\] where $C_1(\beta)$ is a
constant depending on $\beta$ which goes to zero with $\beta$.
\end{Lem}

\pf This is just a computation:

\begin{equation*}
\begin{split}
{d(u_1,u_2) \over t} &\le {d(u_1,y)+d(y,yh_te^{v})+d(yh_te^{v},u_2) \over
t} \text{ 
(by the triangle inequality)} \\
&\le {L\cc(u_1,y)+d(y,yh_te^{v})+L\cc(yh_te^{v},u_2) \over t} \text{
(by L-Lipschitz)}\\
&\le L A \beta +\ri(x,v) + \epsilon + L(C(\beta)) \text{
(by choice of
$\mathcal{B}$ and equation (\ref{spread}))}
\end{split}
\end{equation*}
Let $C_1(\beta)=L(A \beta + C(\beta))$. $\qed$

This end estimate provides the needed precision to make the later
estimates.  The ``well controlled smaller boxes'' mentioned above will
be members of $\mathcal{B}$.  
Other than the end estimates, the other technical point we need is a
measure comparison.  We wish to use the estimate boxes to prove a
relation between $\ri$ and $\rs$ almost everywhere with respect to
$\Ha^Q$.  Thus, we need to know the relative measure of the boxes in
the CC balls.  In light of the discussion concerning the comparison of
Hausdorff and Riemannian measures in section
\ref{CC}, this is fairly straightforward.

\begin{Lem}\label{density}  Fix $T>0$, $\beta >0$ and $v \in
\mathfrak{n}$.  There is a number, $R(\beta,v,T)$, such that $Box(y,tv,t\beta)
\subset B_{cc}(y,tR(\beta,v,T))$ for all $t<T$.  Then,
if $t<T$, the
estimate boxes, $Box(y,tv,\te)$, 
have bounded density variation in CC balls.  i.e. there exist constants
$C_1,C_2$ depending only on $v,\beta$ and $T$ such that
\[ 0<C_1 \le {\Ha^Q(Box(y,tv,t\beta)) \over
\Ha^Q(B_{cc}(y,t(R(\beta,v,T))))} \le C_2 < \infty. \] 
\end{Lem}
\pf  Since $\cc$ is a left invariant metric, the constant
$R(\beta,v,T)$ exists and is independent of the point $y$.  It can
be estimated by equation (\ref{spread}) and the construction of the
boxes.  Since $Q$ is 
the Hausdorff dimension of $(N,\cc)$ and $\Ha^Q(h_t A)=t^q\Ha^Q(A)$, 
$\Ha^Q(B_{cc}(y,tR(\beta,v,T))) = \alpha_Q (tR(\beta,v,T))^Q$ where $\alpha_Q=\Ha^Q(B_{cc}(e^0,1))$.  Thus
we only need estimate $\Ha^Q(Box(y,tv,t\beta))$.  By
construction, $Box(y,tv,t\beta)$ is fibered by paths (radial
geodesic segments)
$F_z=zh_se^v|_{s \in [0,t]}$ where $z \in End(y,tv,t\beta)$.
Now, $\Ha^1(F_z)=\alpha_1t$ (note this is independent of $z$ by the left
invariance of the Hausdorff measure).  Here $\alpha_1$ depends on the
normalization of the Hausdorff measure and $\cc(e^0,e^v)$.  

Let $<\cdot,\cdot>_R$ be a left invariant Riemannian metric on $N$
making the grading orthogonal and matching the sub-Riemannian metric's
inner product on $\mathcal{V}$.  Let $dvol^n$ be the volume form on
$N$ associated to $\ip_R$ and let $dvol^{n-1}$ be the volume form of
$\ip_R$ restricted to a submanifold of one less dimension.  Since
$\Ha^Q$ bounded above and below by a constant multiple of $dvol$, there exists a constants $C_1',C_2'$
such that $0< C_1' \int_{Box(y,tv,t\beta)} dvol \le 
\Ha^Q(Box(y,tv,t\beta))\le  C_2' \int_{Box(y,tv,t\beta)} dvol<
\infty$.  By construction, the 
Riemannian arclength of $F_z$ for $z \in End(y,tv,t\beta)$ is
equal to $t\cc(e^0,e^v)$.  Let $T_yN$ be the (Riemannian) tangent
space to $N$ at 
$y$ and let $R_T \subset \R^{n-1} \subset T_yN$ be the set of vectors such
that $End(y,Tv,T\beta)=ye^{R_T}$.  Now, $R_T \cross [0,T]$ is diffeomorphic to
$Box(y,Tv,T\beta)$ via the map $(w,s) \mapsto ye^we^{sv}$.  Note
that for any $t \le T$ we let $R_t \subset R_T$ be the vectors such that
$ye^{R_t}=End(y,tv,t\beta)$, the same map restricted to $R_t \cross
[0,t]$ is a diffeomorphism onto $Box(y,tv,t\beta)$.  Using Lebesgue
measure on $\R^{n-1}\cross \R$ we see that
\[ \int_{[0,t]}\int_{R_t}d\mathscr{L}^{n-1}d\mathscr{L}^1 =
(2t\beta)^{Q-1}t\cc(e^0,e^v) \]
Using the Jacobian of the diffeomorphism to compute
$\int_{Box(y,tv,t\beta)}dvol$ we get the estimate: 
\[ D_1t(2t\beta)^{Q-1} \le
\int_{Box(y,tv,t\beta)}dvol \le D_2
t(2t\beta)^{Q-1}\]
Where the $D_i$ are constants which take into account the
diffeomorphism's Riemannian distortion of $R_t \cross [0,t]$. 
Note that $D_1>0$ because the Jacobian of the
diffeomorphism can be taken to be bounded below by a positive
constant due to the compactness of $R_T \cross [0,T]$.  Similarly, $D_2$ is finite.  Thus, taking
$C_1=\frac{C_1'D_1(2\beta)^{Q-1}}{C_0(R(\beta,v,T))^Q}$ and
$C_2=\frac{C_2'D_2(2\beta)^{Q-1}}{C_0(R(\beta,v,T))^Q}$, we have the desired
result. $\qed$
\newline

\noindent
{\em Remark:}  Lemma \ref{density} and \ref{Vitali} show that we may refine covers of
the boxes defined above in the same way that we refine covers of
closed balls.  In other words, if $A$ is a measurable set in $N$ and
$\mathscr{B}$ is a cover of $A$ by boxes which is fine at every point
of $A$, then we may find a countable disjoint subcover $\mathscr{B}'$
such that $\Ha^Q(A \minus \cup_{B \in \mathscr{B}'} B)=0$.  This
follows from standard covering arguments (i.e. a countable number of
iterations of the ``greedy''
algorithm used to cover at least a fixed portion of $A$ in one
iteration).  Following Federer, we say that such a collection of boxes
is a {\em Vitali covering realtion}.

\subsection{Differentiability of distances Lipschitz to $\cc$}
The goal of this discussion is to show:
\begin{Thm}\label{KScc}  Let $(N,\mathcal{V},\ip,\cc)$ be a
sub-Riemannian group of nilpotency degree $k+1$.  Suppose 
$d$ is a distance function on $U \cross U$, $U$ an open subset of $N$,
which is L-Lipschitz to $\cc$.  Then, there exists a measurable function
$\Delta_x:N \cross N \ra \R$ defined for all $x \in U$.  Further,
there is a subset
$\mathcal{U}$ of $U$ 
of full measure such that for every $x \in
\mathcal{U}$,
\begin{itemize}
\item For $v \in \mathcal{V}$,
\begin{equation*} {\lim}_{t \ra 0^+} \biggl \{ \bigg|{d(y,yh_te^v) \over t} -
\Delta_x(e^0,e^v) \bigg| \bigg | y \in B_{cc}(x,t)  \biggr \} =0 
\end{equation*}
\item For $v \in \mathcal{V}, \Delta_x(e^0,h_te^v)=t
\Delta_x(e^0,e^v)$.
\item For $v \in \mathcal{V}$,
$\Delta_x(e^0,e^v)=\Delta_x(e^v,e^0)=\Delta_x(e^0,e^{-v})$. 
\item For $v,w, \in \mathfrak{n}$ such that $-v \circledcirc w \in
\mathcal{V}$, $\Delta_x(e^v,e^w)=\Delta_x(e^0,e^{-v}e^w)$.
\end{itemize}
\end{Thm}
\noindent
We will also see that if, for some $x \in \mathcal{U}$,
$(C_xN,\overline{d})$ exists, we can interpret this proposition as a
statement about $\overline{d}$.  One should view the function
$\Delta_x$ as being something like an infinitesimal sub-Finsler
metric and think of the last three claims of the theorem  as giving a limited
homothety, symmetry and left invariance for the function $\Delta_x$.  It would be interesting to find conditions under which this
construction actually converged to a sub-Finsler metric.

\noindent
{\em Remark: }  This theorem leaves open the question of complete
metric differentiaibilty in this setting and other settings.  In a
subsequent paper, the author will address such issues.

We prove theorem \ref{KScc} through several lemmas.  

\begin{Lem}\label{diff1}
If $v \in \mathcal{V}$ and $x$ an approximate continuity
point of $\ri(\cdot,v)$ and $\tau \in \R$, then $\rs(x,\tau v) \le
|\tau| \ri(x,v)$.  In particular, for 
such points $x$, $\ri(x,v)=\rs(x,v)$.
\end{Lem}
\pf  Fix
$\epsilon_0>0$.    We aim to show that for small enough $t>0$ and $y \in
B_{cc}(x,t)$, 
\begin{equation}\label{aim}
{d(y,yh_th_\tau e^v) \over t} \le |\tau|\ri(x,v) + \epsilon_0 
\end{equation}
For simplicity, we will first consider only $\tau>0$.  The nonpositive case
is an easy consequence of the positive case and will be discussed at
the end of the proof.  To prove the positive case, we will consider
special boxes which will easily  
estimate the middle segment of $Box(y,t \tau v,t\beta)$ (for
some appropriate $\beta$) and
still have their endpoint quotients close to the infimum $\ri$.
First, pick $\epsilon, \beta \in (0,\infty)$.  Since $x$ is an
approximate continuity point of $\ri(\cdot,v)$, the set $Z=\{z \st
|\ri(z,v)-\ri(x,v)| < \epsilon\}$ has density in $B_{cc}(x,r)$
approaching $1$ as $r$ approaches $0$.  Considering $y \in
B_{cc}(x,t)$ and $c<1$ such that $B_{cc}(y,ct) \subset B_{cc}(x,t)$,
we observe that the left invariance of $\Ha^Q$ implies that
$\frac{\Ha^Q(B_{cc}(x,t))}{\Ha^Q(B_{cc}(y,ct))}=\frac{1}{c^Q}$ for any
such $y$ and $c$.  Since $Z$ has 
density approaching 1 (with $t$) in $B_{cc}(x,t)$ then $\Ha^Q(Z \cap
B_{cc}(x,t)) \ge (1-\epsilon(t)) \Ha^Q(B_{cc}(x,t))$ where
$\epsilon(t) \ra 0$ as $t \ra 0$.  Thus, since $\Ha^Q(Z \cap
B_{cc}(y,ct)) \ge \Ha^Q(B_{cc}(y,ct))-\epsilon(t)\Ha^Q(B_{cc}(x,t))$,
$\frac{\Ha^Q(Z \cap B_{cc}(y,ct))}{\Ha^Q(B_{cc}(y,ct))} \ge 1- 
\frac{\epsilon(t)}{c^Q}$.  Combining this with the density estimate of
lemma \ref{density} yields that $Z$ has density approaching $1$ in
boxes of the form $Box(y,t\tau v,t\beta)$ with $y \in B_{cc}(x,t)$ as
$t \ra 0^+$.  For the rest of the
proof, fix $t_0$ such that for $t \le t_0$, $y \in B_{cc}(x,t)$,
\[ {\Ha^Q(Z \cap Box(y,t\tau v,t \beta)) \over
\Ha^Q(Box(y,t\tau
v,t \beta))} \ge 1-\epsilon \]

To get estimate (\ref{aim}), we fix a box $B_0=Box(y,t\tau v,t
\beta)$ with $t \le t_0$ and $y \in B_{cc}(x,t)$.  Let $Z_0=Z\cap
B_0$.  Next, let
$\mathcal{B_Z}=\{B \in \mathcal{B}(z,v,\beta, \epsilon) \st B \subset B_0, z
\in Z_0\}$.  Recalling the definition of $\ri$, we see that for $z_0 \in
Z_0$ and if $t$ is sufficiently small (i.e. $t$ such that
$\frac{d(y,yh_te^{v})}{t}$ 
must be close to infimal and $y \in B_{cc}(z_0,t)$), then there are
infinitely many boxes in $\mathcal{B}(z_0,v,\beta, 
\epsilon)$ with height $t\beta$.  Using lemma
\ref{density}, we see that $\mathcal{B_Z}$ is a Vitali covering
relation (see the remark after \ref{density}) of
$Z_0$.  Thus, we may
select from the cover a disjoint 
finite subcover such that $\Ha^Q(\cup_{B \in \mathcal{B}_Z}B) > \Ha^Q(B_0) -
2\epsilon\Ha^Q(B_0)$ which, by abuse of notation, we 
denote this by $\mathcal{B_Z}$ as well.  Consider the picture at this point, we
have fixed a particular box and covered it with smaller boxes, all
contained in the fixed box, such that the end estimates on the small
boxes are approximately infimal.  Next, we will extract a particular curve
from the fixed box and approximate the difference quotient of the
endpoints of its middle segment using the cover. 

Using the same idea as in lemma \ref{density}, we realize $B_0$
as the diffeomorphic image of $R \cross [0,t\tau] \subset \R^{n-1}
\cross \R$.  We can calculate $\Ha^Q(B_0)$, up to bounded error,
by calculating the Lebesgue measure of $R \cross [0,t\tau]$.  Using
this diffeomorphic identification, we denote by $\mathcal{R}_Z$ the
image of the collection 
$\mathcal{B}_Z$ in $R \cross [0,t \tau]$.  The fact that
${\Ha^Q(\cup_{B \in \mathcal{B}_Z} B \cap B_0) \over \Ha^Q(B_0)} \ge
1-2\epsilon$ implies that $\frac{\mathscr{L}^n(\cup_{S \in
\mathcal{R}_Z}S \cap (R \cross [0,t\tau]))}{\mathscr{L}^n(R \cross
[0,t\tau])} \ge 1- 2 \widetilde{C}_0 \epsilon$ where $\widetilde{C}_0$ is a constant depending
on the structure constants of the Lie group and the Jacobian of the
diffeomorphism restricted to the compact set $R \cross [0,t\tau]$.
Thus, using Fubini's  
theorem on $\R^{n-1} \cross \R$, we can conclude that there exists a
point $r_0 \in R$ such that $\frac{\mathscr{L}^1(\cup_{S \in \mathcal{R}_Z}
S \cap r_0 \cross [0,t \tau])}{\mathscr{L}^1(r_0 \cross [0,t\tau])} \ge
1-2 \widetilde{C}_0\epsilon$.  Let $\gamma(t)=y'h_th_\tau e^v$ be the fiber which is
diffeomorphically identified with $r_0 \cross [0,t\tau]$.  Then,
\[{\Ha^1(\cup_{B \in \mathcal{B}_Z} B \cap \gamma) \over
\Ha^1(\gamma)} \ge 1-2C_0\epsilon\]  where $C_0$ is a (potentially) different constant depending on the same data.
Now, we will do end estimates on the smaller boxes and sum them to get
an estimate for ${d(y',y'h_th_\tau e^v) \over t}$.  Let $e_i$,
$i=1,2,3,...,n$ be 
the finite list of points along $\gamma$, listed in order following
the parameter, such that each $e_{2i-1}$ lies in the end of one of the
small boxes and $e_{2i}$ lies in the other.  For $e_i,e_{i+1}$ a pair
which are in the ends of a single 
box, $Box(z_i,{t_i} v,t_i \beta)$, the construction yields the following
estimate.  By lemma \ref{endest}, we have:
\begin{equation}\label{smallbox}
{d(e_i,e_{i+1}) \over t_i} \le C_1(\beta) + \ri(x,v) +  \epsilon 
\end{equation}

\begin{figure} 
\mbox{\epsfig{file=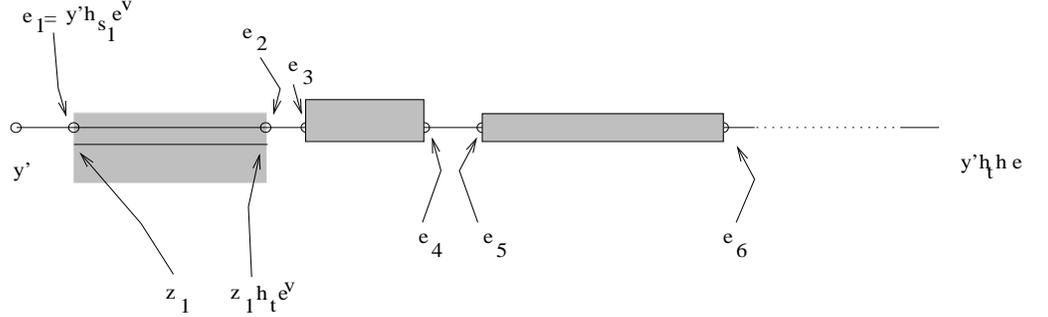,angle=270}}
\caption{A schematic of the covering of the line segment by small boxes.} \label{coverline}
\end{figure}

\noindent
Now using the triangle inequality for $d$, the fact that the smaller boxes cover most of $B_0$ shows that
$d(y',y'h_th_\tau e^v) \le \sum_{i=1}^{\frac{n}{2}}
d(e_{2i-1},e_{2i}) + error$.  Let
$A=\gamma \minus \cup_{B \in \mathcal{B}_Z} B$.  The error is the sum
of the distances 
between the endpoints of the adjacent boxes.  This error is
less that the one dimensional d-Hausdorff measure of $A$,
$\Ha^1_d(A)$ because for any connected set $C$, $\Ha^1_d(C) \ge
diam_d(C)$ (see \cite{Fed} 2.10.12).  Using the Lipschitz property
of $d$, $\Ha^1_d(A) \le L 
\Ha^1(A)$.  Next, we estimate $\Ha^1(A)$.  Since
$\frac{\Ha^1(\cup_{B \in \mathcal{B}_Z} B \cap \gamma)}{\Ha^1(\gamma)}
\ge 1-2C_0\epsilon$, we have that
$\frac{\Ha^1(A)}{\Ha^1(\gamma)}=1-\frac{\Ha^1(\cup_{B \in
\mathcal{B}_Z} B \cap \gamma)}{\Ha^1(\gamma)} \le 1-(1-2C_1\epsilon)$
and using the fact that $\Ha^1_d(A) \le L 
\Ha^1(A)$ and the Lipschitz property of $d$, 
$error \le \Ha^1_d(A) \le L \Ha^1(A) \le 2 C_1\epsilon L \Ha^1(\gamma)
\le 2 C_1 \epsilon \alpha_1 t \tau L \cc(e^0,e^v)$
Letting
$C_2(\epsilon)=2C_1\epsilon \tau L \cc(e^0,e^v)$ we have:
\begin{equation*}
\begin{split}
d(y',y'h_th_\tau e^v) &\le (\sum_{i=1}^{\frac{n}{2}} d(e_{2i-1},e_{2i}))+tC_2(\epsilon)\\
&\le (\sum_{i=1}^{\frac{n}{2}} t_i)(C_1(\beta) + \ri(x,v) +
\epsilon)+tC_2(\epsilon)  
\text{  (by inequality (\ref{smallbox}))}
\end{split}
\end{equation*}
Next, we note that ${\sum_{i=1}^{\frac{n}{2}} t_i \over t}
\le \tau$.  To see this, note that since the curve $\gamma(t)$ is
rectifiable with respect to $\cc$\footnote{This use of rectifiability
is what prevents the same argument from proving differentiability in
nondistributional directions.}, we can pick time indices $s_i$ and
$s_{i+1}$ such that \[\sum_{i=1}^{\frac{n}{2}} t_i = \sum_{i=1}^{\frac{n}{2}}
\frac{\cc(z_i,z_ih_{t_i}e^v)}{\cc(e^0,e^v)} =
\sum_{i=1}^{\frac{n}{2}}\frac{\cc(y'h_{s_i}e^v,y'h_{s_{i+1}}e^v)}{\cc(e^0,e^v)}\]
By left invariance, this equals $\sum_{i=1}^{\frac{n}{2}}
|s_{i+1}-s_i| \le t\tau$.  (See figure \ref{coverline}).
  We use the
estimate and this fact to investigate the 
desired quotient. 
\begin{equation*}
\begin{split}
{d(y',y'h_th_\tau e^v) \over t} &\le |\tau|(C_1(\beta) + \ri(x,v) + 
\epsilon)+C_2(\epsilon) 
\intertext{Using the Lipschitz property as in lemma \ref{endest}, we see:}
{d(y,yh_th_\tau e^v) \over t} &\le |\tau|(C_1(\beta) + \ri(x,v) +
 \epsilon)+ C_2(\epsilon)+C_1(\beta)\\
\intertext{Picking $\beta$ small enough this yields,}
{d(y,yh_th_\tau e^v) \over t} &\le |\tau|\ri(x,v) +
\frac{\epsilon_0}{2} + C_3(\epsilon) 
\end{split}
\end{equation*}
Where $C_3(\epsilon)$ depends only on $\epsilon$.  Thus, picking
$\epsilon$ so that $C_3(\epsilon) < \frac{\epsilon_0}{2}$  allows
us to conclude that ${d(y,yh_th_\tau e^v) \over t}$ $\le |\tau|\ri(x,v) +
{\epsilon_0}$
Taking the $\limsup$ and letting
$\epsilon_0$ tend 
to zero, we have the desired result.  Recall that we had assumed that
$\tau$ was positive at the outset of the proof.  To see that the same
result holds true for negative $\tau$, we make the following
observations.  First, 
$h_{-\tau}e^v=h_\tau e^{-v}$ since $v \in \mathcal{V}$.  Second, for small enough
$t$, the density of $Z$ in $Box(y,-t \tau v,t \alpha)$ is
still greater than $1-\epsilon$.  Therefore, by the symmetry of the
distance $d$, the same argument goes through for $e^{-v}$.
Since the case of $\tau=0$ is trivial, the lemma is valid for all $\tau
\in \R$.  $\qed$ 

\begin{Lem}\label{itsanorm} Retaining the assumptions of theorem
\ref{KScc}, define $\rho(x,\cdot)=\rs(x,\cdot)$.  Then, there exists a
$\mathcal{U}\subset U$ of full measure such that for $x \in
\mathcal{U}$, $\rho(x,\tau v)=|\tau| \rho(x,v)$ for all $v \in
\mathcal{V}$.  We will call $\mathcal{U}$ the set of {\bf metric regular}
points of $d$.  
\end{Lem}
\pf Let $n_0=\dim \mathcal{V}$.  Let $Q$ be a countable subset
of $\mathcal{V}$ so that $e^Q$ is dense in $e^\mathcal{V}$.  Since $\ri$ is measurable function, it
follows that for fixed $v \in \mathcal{V}$, the 
set \[U_v=\{x \st x \text{ is an approximate continuity point of }
\ri(\cdot,v)\}\] has full measure in $U$.  Therefore, $\mathcal{U}=\{x \st x
\in U_v \text{ for all } v \in Q\}=\cap_{v\in Q} U_v$ is a Borel 
set of full measure in $U$.  On this set, we have the desired
properties.  Using lemma \ref{diff1}, we see that for $v\in Q$
and $\tau \in \Q$,
\begin{alignat*}{2}
\ri(x,\tau v) &\le \rs(x,\tau v) &\le |\tau|\ri(x,v) \\
\ri(x,{1 \over \tau} v) &\le \rs(x,{1 \over
\tau} v) &\le \big |{1\over
\tau} \big |\ri(x,v) 
\end{alignat*}
\noindent therefore,
\begin{alignat*}{1}
\ri(x,\tau v) &= |\tau| \ri(x,v) \\
\end{alignat*}

Given the monotonicity in $\tau$ of $\ri$ and $\rs$, this holds for
all $\tau$ in $\R$ and $v \in Q$.  Note in particular
that $\ri(x,-v)=\ri(x,v)$ and therefore we have
symmetry of this type.  To extend these properties to all of $\mathcal{V}$ we
check that $\rho(x,\cdot)$ is continuous.  Fix $v \in \mathcal{V}$
and consider a sequence $v_i \in Q$ which converges to $v$.  Then,
\begin{equation*}
\begin{split}
\rho(x,v_i)-\rho(x,v) &=\overline{\lim}_{t \ra 0^+}
\bigg \{\frac{d(x',x'h_te^{v_i})}{t} \bigg | x' \in B_{cc}(x,t) \bigg
\} - \\ & \qquad \overline{\lim}_{t \ra 0^+}
\bigg \{\frac{d(x'',x''h_te^v)}{t} \bigg | x''
\in B_{cc}(x,t)\bigg \}\\
&\le \overline{\lim}_{t \ra 0^+}
\bigg \{\frac{d(x',x'h_te^{v_i}) - d(x',x'h_te^v)}{t} \bigg | x'
\in B_{cc}(x,t)\bigg \}\\
& \le \overline{\lim}_{t \ra 0^+} \bigg
\{\frac{d(x'h_te^{v_i},x'h_te^v)}{t}\bigg | x' \in B_{cc}(x,t)\bigg \}
\\
&\le L \overline{\lim}_{t \ra 0^+} \cc(e^{v_i},e^v)
\intertext{Which tends to zero as $i \ra \infty$.  A similar argument
shows that $\lim_{i \ra \infty} \rho(x,v_i) - \rho(x,v) \ge 0$ as well.}
\end{split}
\end{equation*}
$\qed$
\begin{Lem}\label{lirho} Suppose $x_0 \in N$ is a metric regular point.
Then, if $v,w \in 
\mathfrak{n}$ such that $-v \circledcirc w \in \mathcal{V}$ then  
\[\overline{\lim}_{t \ra 0^+} \biggl \{{d(yh_t e^v,yh_t e^w) \over t}
\bigg | y \in B_{cc}(x_0,t) \biggr \} = \rho(x_0,-v \circledcirc w) \]
\end{Lem}
\pf This fact follows from the covering argument in lemma
\ref{diff1}.  Fix $\epsilon_0>0$.  Using the same notation, picking $\beta$
sufficiently small, the density of $Z_0$ in
$Box(x_0h_te^v,t(-v \circledcirc w),\beta)$ is 
bigger than $1-\epsilon$.  Now, we can run the same covering argument
to produce the estimate ${d(x_0h_te^v,x_0h_te^w) \over t} \le \ri(x_0,-v
\circledcirc w) + \epsilon_0$.$\qed$ \newline

\noindent
{\em Proof of theorem \ref{KScc}:}  Let
\[\Delta_x(e^{v},e^{w})=\overline{\lim}_{t \ra 0^+} \biggl \{
\frac{d(yh_te^v,yh_te^w)}{t} \bigg | y \in B_{cc}(x,t) \biggr \}\]
This clearly satisfies the first part of the claim with $\mathcal{U}$
being the set of metric regular points of $d$.  We now prove the other items in
reverse order.  By lemma \ref{lirho}, we see that
$\Delta_x(e^v,e^w)=\Delta_x(e^0,e^{-v \circledcirc w})$ for $-v \circledcirc w
\in \mathcal{V}$ proving the fourth part of the theorem.
  Lemma \ref{itsanorm} implies that
$\Delta_x(e^0,e^v)=\Delta_x(e^0,e^{-v})$ and the symmetry of the
distance function $d$ immediately shows that $\Delta_x(e^0,e^v)=\Delta_x(e^v,e^0)$, proving the third part of the claim.  Lemma
\ref{itsanorm} also shows that for $v \in \mathcal{V}$,
$\Delta_x(e^0,h_te^v)=t\Delta(e^0,e^v)$ proving the second part of the
theorem.  $\qed$
\section{Tangent cones}\label{Tangcones}
The differentiability result of the last section allows us to compare
the local structures of a \newline sub-Riemannian group and some metric
space to which it is biLipschitz equivalent.  Presuming the existence
of tangent cones to image points of $F$ in $X$, we would like construct
a ``derivative'' mapping between the tangent cone at a metric regular
point in $N$ 
to the tangent cone at its $F$ image in $X$.  To do this, we will
first review the tangent cone construction and apply it to this
setting, interpreting the results of theorem \ref{KScc}.  The goal
of this section is to prove the following:
\begin{Pro}\label{conemap} Let $(N,\mathcal{V},\ip,\cc)$ be a
sub-Riemannian group, $(X,d_X)$ be a metric space and $F:N \ra X$ be
an L-biLipschitz map between them.  Letiing $d=F^*d_X$ we let
$\Delta_X$ be the function constructed in Theorem \ref{KScc}.  Suppose
there exists an $x \in X$ 
in the image of the metric regular points of $F$ such that a tangent cone,
$(C_xX,\overline{d}_X)$, exists.  Then,
\begin{itemize}
\item There exists an L-biLipschitz map
$\overline{F}:(C_{F^{-1}(x)}N,\overline{d}_{cc}) \ra
(C_xX,\overline{d}_X)$ such that for $n^{-1}m \in C_{F^{-1}(x)}N$ that is
in the exponential image of $\mathcal{V}$,
$\overline{d}_{X}(\overline{F}(n),\overline{F}(m))=\Delta_{F^{-1}(x)}(n,m)$.
\item $\overline{F}$ takes radial geodesics to geodesic lines.
\end{itemize}
\end{Pro} 
We begin by considering the structure of the tangent cone to $(N,\cc)$
at some point.  By proposition \ref{Mitchell}, we
know that the tangent cone at any point is unique and isometric to
$(N,\cc)$.  We will consider the tangent cone at $e^0$ for
simplicity.  Given a complete metric space $(X,d)$, one usually
constructs the tangent cone at $x$ as the Gromov-Hausdorff limit of the
sequence of dilated spaces $(X,x,\lambda_id)$ for some sequence
$\{\lambda_i\}$ tending towards infinity as $i \ra \infty$.  For
$(N,\cc)$, one may construct the dilated spaces using the homothety.
Thus one considers the convergence of $(N,h_{\lambda_i}^*\cc)$.  Now,
since $(N,\cc)$ is isometric to $(N,h_{\lambda_i}^*\cc)$ via the
isometry $h_{\lambda_i}$, the sequence is, in fact, the constant
sequence and thus converges in the Gromov-Hausdorff topology for any
choice of sequence $\{\lambda_i\}$.  (See \cite{GLP} for an explanation for the
convergence of pointed spaces.  The reader may also want to consult
\cite{GHreview} which provides a somewhat shorter  explanation for compact
spaces).  
Before continuing, we reiterate the assumptions about $(X,d_X)$:
first, there exists $x \in Im(F) \subset X$ and the tangent cone to X
at x exists, i.e. there exists a sequence $\{\lambda_i\}$ such that
$(C_xX,\overline{d}_X)= \text{GH-lim}_{i \ra \infty}(X,x,\lambda_i
d_X)$ (the Gromov-Hausdorff limit of $(X,x,\lambda_id_X)$).  To
construct the map $\overline{F}$, we first need to review 
the concept and usage of ultrafilters and ultralimits.
\begin{Def}  A nonprincipal ultrafilter is a
finitely additive probability measure $\omega$ on subsets of
$\mathbb{N}$ such that 
\begin{itemize}
\item $\omega$ is zero on any finite subset of  $\mathbb{N}$.
\item For every subset $S \subset \mathbb{N}$, $\omega(S)$ is either
zero or one.
\end{itemize}
\end{Def}
It is a consequence of the axiom of choice that nonprincipal
ultrafilters exist (see \cite{ultraref} chapter 10, theorem 7.3).  
We can use an ultrafilter to choose a distinguished ``convergent'' point for a
sequence.  As an illustrative example, let $\R$ be the line with its
usual metric  and let $\{x_i\}$ be a 
bounded sequence in $\R$.  Think of $\{x_i\}$ as the image of a map
$s:\mathbb{N} \ra \R$.  Then there  is a unique element $x_\omega \in
X$ ($x_\omega = \wlim \text{ }s$) such that  for any neighborhood $U$
of $x_\omega$, 
$\omega(s^{-1}(U))=1$.  In this way, the nonprincipal ultrafilter
picks a limit point of the sequence.  Next we point out a useful fact about
ultralimits of 
sequences in $\R$.  If $\{a_i\}$ and
$\{b_i\}$ are bounded sequences and $a_i 
\le b_i$ then $a_\omega=\wlim \text{ }a_i \le b_\omega=\wlim \text{
}b_i$.  Indeed, if 
$b_\omega < a_\omega$ then we can find open neighborhoods $U_a$ and
$U_b$ of $a_\omega$ and $b_\omega$ which separate the ultralimits.
Since $\omega(s^{-1}(U_a))=1=\omega(s^{-1}(U_b))$, if the inverse
images are disjoint, this violates the additivity of the ultrafilter
while if they intersect they must intersect in an infinite set by the
additivity of the ultrafilter.   Thus, for indices in this intersection,
$\lim b_i < \lim a_i$, contradicting the assumption that $a_i \le b_i$ for
all $i$.  

Let $\omega$ be a nonprincipal ultrafilter and define a limit metric
space $(C_x^\omega X,x,d_X^\omega)$ as follows.  Let
$\overline{X}=\{\{x_i\} | x_i \in (X,x,\lambda_id_X),
\lambda_id_X(x_i,x) < \infty \}$. Then $C_X^\omega X= \overline{X} /
\sim$ where $\{x_i\} \sim \{y_i\}$ if \newline $\wlim_{i \ra \infty} \lambda_i
d_X(x_i,y_i)=0$ (note that here $d_X(x_i,y_i)$ is a bounded map from
$\mathbb{N}$ to $\R$) and $d^\omega_X(\{x_i\},\{y_i\})=\wlim_{i \ra \infty}
\lambda_i d_X(x_i,y_i)$.  Note that the choice of ultrafilter forces
convergence of the dilated spaces and allows us to group the many
 choices we otherwise have to make into the choice of the ultrafilter.  Also note that our
assumption that the Gromov-Hausdorff limit of $(X,x,\lambda_id_X)$
exists allows us to find an isometry between $C_xX$
and $C_x^\omega X$.  For a proof of this fact, as well as a more
general discussion of ultrafilters and ultralimits, see
\cite{KleinerLeeb} section 2.4.2 and, in particular, lemma 2.4.3.  

Using left invariance, we may assume without loss of generality that
$F^{-1}(x)=e^0$.  We 
define $\overline{F}:(C_{e^0}N,e^0,\overline{d}_{cc}) \ra
(C_x^\omega X,x,d^\omega_X)$ by
$\overline{F}(n)=\{F(h_{\lambda_i^{-1}}n)\}$.  Note that
$\overline{F}$ is L-biLipschitz if $F$ is:
\begin{equation*}
\begin{alignat*}{2}
\intertext{ By definition,}
d^\omega_X(\overline{F}(n),\overline{F}(m))&=\wlim_{i \ra \infty}
\lambda_i d_X(F(h_{\lambda_i^{-1}}n), 
F(h_{\lambda_i^{-1}}m))\\
\intertext{Using the biLipschitz property,}
\wlim_{i \ra \infty} \frac{1}{L}\cc(n,m)&  \le \wlim_{i \ra \infty}
\lambda_i d_X(F(h_{\lambda_i^{-1}}n), 
F(h_{\lambda_i^{-1}}m)) & \le \wlim_{i \ra \infty} L \cc(n,m)\\
\intertext{Upon taking ultralimits, we have:}
\end{alignat*}
\end{equation*}
\begin{gather*}
\frac{1}{L}\overline{d}_{cc}(n,m) \le d^\omega_X(\overline{F}(n),
\overline{F}(m)) \le L \overline{d}_{cc}(n,m) 
\end{gather*}

As mentioned above, lemma 2.4.3 in \cite{KleinerLeeb} shows that
$(C_xX,\overline{d}_X)$ is isometric to $(C_x^\omega,d_X^\omega)$ and
therefore we may view $\overline{F}$ as an L-biLipschitz map between
$(C_{e^0}N,\overline{d}_{cc})$ and $(C_xX,\overline{d}_X)$.
In addition, by construction,
$\overline{d}_X(\overline{F}(n),\overline{F}(m)) = 
\wlim_{i \ra \infty} \lambda_i
d_X(F(h_{\lambda_i^{-1}}n),F(h_{\lambda_i^{-1}}m))=\lim_{i \ra \infty}
\frac{d(h_{\lambda_i^{-1}}n,h_{\lambda_i^{-1}}m)}{\lambda_i}$ which by Theorem
\ref{KScc} applied to $F^*d_X$ converges to $\Delta_x(n,m)$ when
$n^{-1}m$ is in the 
exponential image of $\mathcal{V}$.  Thus we have proved the
first part of proposition \ref{conemap}.  For the second part,
consider a radial 
geodesic, $h_te^v$ in $N$ where $v \in \mathcal{V}$.  Using the
definition of $\overline{F}$, let $\gamma(t)=\overline{F}(h_te^v)$.
We claim $\gamma(t)$ minimizes between any two time indices $t_1 <
t_2$.  Assume for now that $t_1 >0$.  Indeed,
\begin{equation*}
\begin{split}
\overline{d}_X(\gamma(t_1),\gamma(t_2))=\overline{d}^\omega_X(\overline{F}(h_{t_1}e^v),\overline{F}(h_{t_2}e^v))
&= \Delta_x(h_{t_1}e^v,h_{t_2}e^v) =
\Delta_x(e^0,(h_{t_1}e^v)^{-1}(h_{t_2}e^v)) \\
&= \Delta_x(e^0,h_{t_1}e^{-v}h_{t_2}e^v) =
\Delta_x(e^0,h_{t_2-t_1}e^v)\\
&= |t_2 -t_1| \Delta_x(e^0,e^v)
\end{split}
\end{equation*}
Here we used the left invariance of $\Delta_x$ and lemma \ref{ht}.
For $t_1<0$, we simply recall our convention regarding
$h_tn$ for negative values of $t$.  When $t_1,t_2\ge 0$ we observe
that
\begin{equation*}
\begin{split}
\overline{d}_X(\gamma(-t_1),\gamma(t_2))=\overline{d}^\omega_X(\overline{F}(h_{-t_1}e^v),\overline{F}(h_{t_2}e^v))
&= \Delta_x(h_{t_1}e^{-v},h_{t_2}e^v) =
\Delta_x(e^0,(h_{t_1}e^{-v})^{-1}(h_{t_2}e^v)) \\
&= \Delta_x(e^0,h_{t_1}e^{v}h_{t_2}e^v) =
\Delta_x(e^0,h_{t_2+t_1}e^v)\\
&= |t_2 +t_1| \Delta_x(e^0,e^v)
\end{split}
\end{equation*}
Thus, $\gamma(t)$ minimizes
between any two time indices.  In addition, by left invariance, all
radial geodesics are taken to geodesics under the map $\overline{F}$.

The reader may wonder why such trouble is taken to use an ultrafilter
construction to produce convergence.  In fact, if $X$ is a separable
and locally compact metric space, then one can use more standard
(e.g. Arzela-Ascoli type) arguments to produce uniform convergence on
compact sets.  However, since asymptotic cones are not always locally
compact, we use the ultrafilter construction to resolve convergence
issues.  
\section{Generalized notions of curvature}\label{Curv}
\subsection{Curvature bounded above}
Again, we follow the notation and definitions of \cite{KleinerLeeb}.
To determine 
curvature bounds of metric spaces, we will need to 
compare to model spaces.  For $\kappa \in \R$, let $(M^2_\kappa,d_\kappa)$
\index{$M^2_\kappa$} be the two dimensional constant curvature
$\kappa$ model space (i.e. $\R^2, S^2_\kappa$, or
$\mathbb{H}^2_\kappa$).  Denote by 
$diam(\kappa)$ the diameter of the model space $M^2_\kappa$.
Let $(X,d)$ be a complete metric space.  Given a triangle, $\triangle$, in X
with minimizing geodesic sides $s_1, s_2, s_3$, we call the
triangle $\triangle_\kappa \subset M^2_\kappa$ with the same side
lengths as $\triangle$ a comparison triangle for $\triangle$.  
\begin{Def}
$(X,d)$ is a $CAT_\kappa$\index{$CAT_\kappa$} space if
\begin{itemize}
\item For every pair of points $x,y \in X$ with $d(x,y) <
diam(\kappa)$, there is a geodesic segment joining $x$ to $y$.
\item Let $\triangle$ be geodesic triangle in $X$ with sides
$s_1,s_2,s_3$ of lengths $l_1,l_2,l_3$ such that $l_1+l_2+l_3 <
2diam(\kappa)$.  For any two points $x,y \in \triangle$, $d(x,y) \le
d_\kappa(x_\kappa,y_\kappa)$ where $x_\kappa$ and $y_\kappa$ are the
corresponding points on $\triangle_\kappa$.  
\end{itemize}
\index{$CAT_\kappa$}
\end{Def}

Another way of saying the second condition in the definition is that
the triangles in $X$ are thinner than their counterparts in
$M^2_\kappa$.  As mentioned in the introduction, we sometimes will
consider metric spaces which have this triangle property locally.
Hence, we say a metric space, $X$, is a $CBA_\kappa$ metric space if for
each point $x \in X$, there exists an $r$ such that the closure of the
ball of radius $r$ about $x$ is a $CAT_\kappa$ metric space in its own
right.  In the next several lemmas,
we outline the properties of $CAT_\kappa$ 
spaces that will be of use below.  A nicely self-contained
discussion of these facts can be found in \cite{Ballmann:SNPC}.  

\begin{Thm}[Hadamard-Cartan Theorem for ${CAT_0}$ spaces] Let $X$
be a simply connected complete \\$CAT_0$ metric space.  Then, for any
two points $x,y \in X$, there exists a unique geodesic connecting $x$ to $y$.
\end{Thm}

\begin{Lem}\label{CATgrowth}  Let $(X,d)$ be a $CAT_0$ space and
suppose $\gamma_1(t)$ and $\gamma_2(t)$ are unit speed geodesic rays
in $X$.  Then the function $d(\gamma_1(t),\gamma_2(t))$ either grows
linearly in $t$ or is bounded.
\end{Lem}
\pf The key fact is that in a $CAT_0$ space, the distance function is
convex (for a proof, see \cite{Ballmann:SNPC} proposition 5.4).  If
the function $f(t)=d(\gamma_1(t),\gamma_2(t))$ is linear or bounded, we are
done.   Assuming that
the function is not linear but is unbounded, we can find an increasing sequence
$\{t_i\}$ such that $f''(t_i)>0$ for all $i$.  One can now easily show that
$f(t)$ must grow at least linearly in this case.  But, $f(t)$ is
bounded above by $2t+d(\gamma_1(0),\gamma_2(0))$ simply by measuring 
along the two geodesics.  Thus, $f(t)$ grows linearly or is bounded. $\qed$

For the sake of completeness, we quote some results concerning the two
coning operations we use in this argument.

\begin{Lem}\label{CATcone}
Let $(X,d)$ be a $CAT_\kappa$ space.
\begin{enumerate}
\item For all $x \in X$, $C_xX$ exists and is a $CAT_0$ space.
\item If $\kappa =0$ and $\omega$ is a nonprincipal ultrafilter, then
$C^\omega_\infty X$ is a $CAT_0$ space.  
\end{enumerate}
\end{Lem}
\noindent
{\em Proof:} For 1) This is proved in \cite{Nik:tangcone}, 2) see
\cite{KleinerLeeb} or \cite{Ballmann:SNPC}.

Also note that the lemma implies that the tangent cone to a
$CBA_\kappa$ space is also a $CAT_0$ space.

\subsection{Curvature bounded below}
Curvature bounded below for metric spaces is defined in much
the same way as for $CAT_\kappa$ spaces.  A standard resource for this
material is \cite{BGP}.  The definition and first uses of these spaces
were due to A. D.
Alexandrov and are often referred to as ``Alexandrov spaces''.  Using the notation of the
previous section, we make the following definitions.
\begin{Def}  Let $(X,d)$ be a complete metric space.  Then,
$(X,d)$ is a $CBB_\kappa$ space if
\begin{itemize}
\item $(X,d)$ is a locally compact geodesic space.
\item Let $\triangle$ be geodesic triangle in $X$ with sides
$s_1,s_2,s_3$ of lengths $l_1,l_2,l_3$ such that $l_1+l_2+l_3 <
2diam(\kappa)$.  For any two points $x,y \in \triangle$, $d(x,y) \ge
d_\kappa(x_\kappa,y_\kappa)$ where $x_\kappa$ and $y_\kappa$ are the
corresponding points on $\triangle_\kappa$.  
\end{itemize}
\end{Def}
As with $CAT_0$ spaces, there is a natural picture, $CBB_0$ triangles
are ``fatter'' than Euclidean triangles.  As mentioned in the
introduction, by the usual arguments using Toponogov's theorem the
localized version of this definition curvature is equivalent to the
global and, for the purposes of this exposition, do not distinguish
between the two.     
We end this section with some useful theorems concerning the geometry of
$CBB_\kappa$ spaces.
\begin{Lem}\label{CBTcone}  Let $(X,d)$ be a $CBB_\kappa$ space.
\begin{enumerate}
\item Any tangent cone to $X$ at $x$ is a $CBB_0$ space.  
\item If $\kappa=0$ and $\omega$ is a nonprincipal ultrafilter then
$C_\infty^\omega X$ is a $CBB_0$ space.
\end{enumerate}
\end{Lem}
\pf The first claim is discussed in \cite{BGP}.  The discussion in
section 2.4.2 in \cite{KleinerLeeb} reduces the second claim to the
claim that the tangent cone to a $CBB_0$ space is $CBB_0$ as well.
The latter fact is also discussed in \cite{BGP}.  

\begin{Thm}[Grove-Petersen \cite{GP1}]\label{GP}  Let $X$ be a locally compact,
complete, noncompact $CBB_0$ metric space.  If $X$ contains a line (i.e.
an unbounded geodesic which minimizes for all time) then $X$ splits
isometrically as $E \cross \R$ where $E$ is some totally convex subset
of $X$ and $\R$ is the real line with its usual metric.
\end{Thm}

For our purposes, we use this theorem to prove an analogous result to
lemma \ref{CATgrowth}.
\begin{Lem} \label{CBTgrowth}  Let $(X,d)$ be a $CBB_0$ space and
let $\gamma_1(t)$ and $\gamma_2(t)$ be two lines in $X$.  Then,
$\gamma_1(t)$ and $\gamma_2(t)$ are either parallel or diverge
linearly.
\end{Lem}
\pf  Since these are lines, theorem \ref{GP} says that $X$ splits
isometrically as $X_0 \cross \R$ where $\R$ is the geodesic
$\gamma_1$.  Now, $\gamma_2(t)$ is a geodesic line in the product and
hence its projection to each factor must be a geodesic line as well.
Therefore, the projection of $\gamma_2(t)$ in $X_0$ is either a point
or a line.  If it is a point, then $\gamma_1$ and $\gamma_2$ are parallel.
If its projection in $X_0$ is not a point, then since it is a line,
this implies that $d(\gamma_1(t),\gamma_2(t))$ must be linear in $t$. $\qed$ 
\section{Applications of proposition \ref{conemap}}\label{end}
The purpose of this section is to use the differentiability result
(albeit weakly) to examine the possible quasi-isometric embeddings of
$N_0$, 
a nonabelian simply connected connected nilpotent Lie group equipped
with a left invariant Riemannian 
metric.  As discussed in the introduction, we wish to prove that
quasi-isometric embeddings of such a Lie group into $CAT_0$ or $CBB_0$
cannot exist.  After taking asymptotic cones, denoting by $N$ the
resulting nonabelian simply connected connected graded nilpotent Lie
group equipped with a 
left-invariant sub-Riemannian metric,  we reduce to proving the
following theorem.\newline

\noindent
{\bf Theorem C.} {\em Let $G$ be a connectes simply connected
  grapded nilpotent Lie group equipped with a left invariant
  Carnot-Carath\`eodory metric and $U \subset G$ be an open set.  Then
  $U$ does not admit a biLipschitz embedding into any Alexandrov
  metric space with curvature bounded above (resp. bounded below).  }

This follows quickly from the next lemma which estimates the growth
rate between radial geodesics.

\begin{Lem}\label{diverge}  Suppose $f: (N,\mathcal{V},\ip,\cc) \ra
(X,d)$ is an L-biLipschitz map such that radial geodesics in $N$ are
mapped to $d$-geodesic lines in $X$.  Then, there exists a pair of
geodesics $\gamma_1(t)$ and
$\gamma_2(t)$ in $X$ which are images of radial geodesics and
constants $C_1,C_2 >0$ and $0 < \alpha < \beta <1$ such that, 
\[ C_1 |t|^\alpha \le d(\gamma_1(t),\gamma_2(t)) \le C_2 |t|^\beta \]
for $|t| > 1$.
\end{Lem}
\pf This is an exercise in using the Campbell-Baker-Hausdorff
formula.  Suppose $\gamma_1(t)=f(h_te^v)$ and
$\gamma_2(t)=f(e^wh_te^v)$ for $v,w \in \mathcal{V}$ such that $[v,w]
\neq 0$.  Using the 
biLipschitz property, we reduce immediately to consider the rate of
growth of the function $\cc(h_te^v,e^wh_te^v)$ for $|t|>1$.
\begin{equation*}
\begin{alignat*}{2}
\cc(h_te^v, e^{w}h_te^v) &= \cc(e^0,h_te^{-v}e^wh_te^v) \\
&= \cc(e^0,e^{w+t[w,v]+H(w,v,t)})\\
\intertext{Where $H(v,w,t)$ are the higher order bracket terms in the
Campbell-Baker-Hausdorff expansion.  Note that the order $n$ bracket
terms in the expansion each 
have a coefficient which of the form $Bt^m$ where $B$ is a constant
from the Campbell-Baker-Hausdorff expansion and $m < n$.
Therefore, for $0< b <1$ this is equal to,}
& \cc(e^0, h_{{|t|^{b}}}e^{\frac{1}{|t|^{b}}w+ \frac{t}{|t|^{2b}}[w,v]+
H'(w,v,t)})\\ 
&= |t|^b \cc(e^0,e^{\frac{1}{|t|^{b}}w+ \frac{t}{|t|^{2b}}[w,v]+
H'(w,v,t)})\\ 
\intertext{Again, $H'(w,v,t)$ are higher order bracket terms but now,
the terms or order $n$ each have a coefficient 
of the form $\pm Bt^{m-{n b}}$, where $B$ a constant.  Now, picking $0 <
\alpha < \beta <1$ such that $m-n\alpha >0$ for at least one value of
$n$ and $m-n \beta \le 0$ for all values of $n$, we have that:}
C_1 |t|^\alpha &\le \cc(h_te^v,e^wh_te^v) \le C_2 |t|^\beta\\
\end{alignat*}
\end{equation*}
Where $C_1$ is the minimum of the function $t \mapsto
\cc(e^0,e^{\frac{1}{|t|^{\alpha}}w+
\frac{t}{|t|^{2\alpha}}[w,v]+H'(w,v,t)})$ on $(-\infty,-1] \cup
[1,\infty)$ and $C_2$ is the maximum of the function $t \mapsto
\cc(e^0,e^{\frac{1}{|t|^{\beta}}w+ \frac{t}{|t|^{2\beta}}[w,v]+
H'(w,v,t)})$ on the same domain.  $C_1$ exists because the function in question is nonnegative and
goes to infinity  as $t \ra \pm \infty$ while $C_2$ exists because the
function in question is nonegative and goes to zero as $t \ra \pm
\infty$.  Thus, using the 
L-biLipschitz property of the map, we have the desired result with the
constants $C_i$ adjusted by either $L$ or $\frac{1}{L}$. $\qed$

\noindent
{\em Proof of Theorem C:}  Suppose $f:N \ra X$ is the
L-biLipschitz map we wish to investigate and $X$ is either $CBA_\kappa$ or
$CBB_\kappa$.  Then, proposition \ref{conemap} lets us reduce to an
examination of a biLipschitz map $\overline{F}:(N,\cc) \ra
(C_xX,\overline{d})$ where $(N,\cc)$ and $(C_xX,\overline{d})$ are the
appropriate tangent cones and $\overline{F}$ takes radial geodesics to
geodesic lines.  Lemma \ref{CATcone} or lemma \ref{CBTcone} implies that
the tangent cone $C_xX$ is either $CAT_0$ or $CBB_0$.  Now we satisfy
the assumptions of lemma \ref{diverge} and hence have a pair of
geodesics with the specified divergence rate.  However, this
contradicts lemma \ref{CATgrowth} or lemma \ref{CBTgrowth}.$\qed$

\bibliographystyle{alpha}

\end{document}